\documentclass{article}

\usepackage{graphicx}

\newtheorem{theorem}{Theorem}
\newtheorem{proposition}[theorem]{Proposition}
\newtheorem{remark}[theorem]{Remark}
\newtheorem{definition}[theorem]{Definition}

\newtheorem{corollary}[theorem]{Corollary}

\newcommand{\R}{{\rm I}\!{\rm R}}

\begin{document}

\title{On the well-posedness of a mathematical model for Lithium-ion batteries\thanks{Accepted for publication in  {\em Applied Mathematical Modelling}}}

\author{Angel Manuel Ramos \\
\ \\
Instituto de Matem\'atica Interdisciplinar \\ Departamento de Matem\'atica Aplicada \\
Facultad de Ciencias Matem\'aticas \\ Universidad Complutense de Madrid \\
Plaza de Ciencias 3, \\
Madrid, 28040, Spain \\
angel@mat.ucm.es}

\date{May 29th, 2015\thanks{A previous version of this work can be seen here: http://www.mat.ucm.es/deptos/ma/prepublicaciones/2013/2013-08p.pdf}}

\maketitle

\begin{abstract}
In this article we discuss the well-posedness of a mathematical model that is used in the literature for the simulation of Lithium-ion (Li-ion) batteries. First, a mathematical
model based on a macro-homogeneous approach is presented, following previous works.
Then it is showed, from a physical and a mathematical point of view, that a boundary condition widely used in the literature is not correct. Although these errors could be just sign typos (that can be explained
as carelessness over d/d$x$ versus d/d$n$, with $n$ the outward unit vector) and
authors using this model probably use the correct boundary
condition when they solve it in order to do simulations, readers should be aware of the right choice.
Therefore, the deduction of the correct boundary condition and a mathematical study of the well-posedness of the corresponding problem is carried out here.
\end{abstract}

\noindent
{\bf Keywords:} Lithium-ion batteries; mathematical model; partial differential equations,
boundary conditions; uniqueness of solution; existence of solution

\section{Introduction}
Lithium-ion (Li-ion) batteries have become very popular in the last years as a source of energy in multiple portable electronic devices.
A mathematical model showing the key factors of the battery operation can be very helpful for the design and optimization of new models and also for the real
time control of its performance.

Based on a macro-homogeneous approach developed by Newman (see \cite{1973Ne}) several mathematical models have been developed for these
purposes (see \cite{1993DoFuNe,1994FuDoNe,2002GoEtAl,2005FaPl,2006SmWa,2006-b-SmWa,2007SmRaWa,2010ChEtAl,2011GuSiWh,2012KiSmIrPe}) which include the main physics present in
charge/discharge processes. A fully mathematical model is presented in Sec.~\ref{themodel}, including a system of boundary value problems for the conservation of Lithium and conservation of charge
in the solid and electrolyte phases, together with an initial value problem for the conservation of energy. In the literature (see \cite{2010ChEtAl}, \cite{2011GuSiWh}, \cite{2012PuMc}, \cite{2007SmRaWa}, \cite{2006SmWa}, \cite{2006-b-SmWa}) one can find numerical computations of the model (or simplifications of it), with parameters corresponding to actual devices, that help to highlight the structure of this highly coupled model and show the relevance to the applications.

Over the last years, several authors have written articles in different journals including a boundary condition that is not correct (see Remark \ref{rmswbc}). Although the authors probably use the correct boundary
condition when solving the model in order to do simulations, the reader should be aware of the right choice. In this article it is shown why that condition is
not only physically incorrect (see Sec.~\ref{cbcs}) but also mathematically incorrect, since it is proved (see Sec.~\ref{eusol}, Remark \ref{nosolu}) that the corresponding system of boundary value
problems modelling the conservation of charge does not have any
solution (therefore the system is not well-posed).
Some results regarding the uniqueness and existence of solution (with the correct boundary condition) for a simplified version of the model are shown in Sec.~\ref{eusol}.

\section{Mathematical model} \label{themodel}

\subsection{Generalities}

A typical Li-ion battery cell has three regions: A porous negative electrode, a porous positive electrode and an electron-blocking separator. Furthermore, the cell contains an electrolyte, which
is a concentrated solution containing charge species that can move all along the cell in response to an electrochemical potential gradient.

The negative electrode is an intercalated lithium compound usually made from Carbon (typically graphite), with Li$_y$C$_6$ active material.
Here $y\in [0,1]$ is the stoichiometry value of the material, which
changes during charge and discharge. For instance, during discharge, lithium ions inside of solid Li$_y$C$_6$ particles diffuse to the surface where they react and transfer from the
solid phase into the electrolyte phase (see \cite{2006SmWa,2006-b-SmWa}). During charge they follow the opposite way.
Theoretically, in a fully charged Li-ion battery this compound is
LiC$_6$ (lithiated graphite; Li saturation; $y=1$), in a semi charged/discharged battery it is Li$_y$C$_6$, with $y\in (0,1)$ and in a fully discharged battery it is just carbon
(Li depletion; $y=0$). In practical operating cases $y$ never attains the extreme values $y=0$ or $y=1$. The corresponding negative reaction
equation is the following:
$$
y\mbox{Li}^+ + y\mbox{e}^-+6\mbox{C} \longleftrightarrow \mbox{Li}_y\mbox{C}_6 \ \ \ \ (0\leq y \leq 1).
$$

The positive electrode is usually a metal oxide or a blend of multiple metal oxides (see \cite{2010ChEtAl}) such as lithium cobalt oxide (Li$_{1-y}$CoO$_2$), lithium iron phosphate
(Li$_{1-y}$FePO$_4$), or lithium manganese oxide (Li$_{1-y}$Mn$_2$O$_4$), with $y\in [0,1]$. For instance, during discharge the positively charged ions travel via diffusion and
migration through the electrolyte solution to the positive electrode where they react and insert into solid metal oxide particles (see, e.g., \cite{2006SmWa,2006-b-SmWa}).
During charge they follow the opposite
way. Theoretically, in a fully charged Li-ion battery $y=1$ (Li depletion), in a semi
charged/discharged battery $y\in (0,1)$ and in a fully discharged battery $y=0$ (Li saturation). Again, in practical operating cases $y$ never attains the extreme values $y=0$ or $y=1$.
The corresponding positive reaction equation for the examples showed above are the following:
$$
\mbox{LiCoO}_2 \longleftrightarrow \mbox{Li}_{1-y}\mbox{CoO}_2 + y\mbox{Li}^+ + y\mbox{e}^- \ \ \ \ (0\leq y \leq 1),
$$
$$
\mbox{LiFePO}_4 \longleftrightarrow \mbox{Li}_{1-y}\mbox{FePO}_4 + y\mbox{Li}^+ + y\mbox{e}^- \ \ \ \ (0\leq y \leq 1),
$$
$$
\mbox{LiMn}_2\mbox{O}_4 \longleftrightarrow \mbox{Li}_{1-y}\mbox{Mn}_2\mbox{O}_4 + y\mbox{Li}^+ + y\mbox{e}^- \ \ \ \ (0\leq y \leq 1).
$$
Therefore, considering both, negative and positive reaction, the corresponding total reaction equations are the following:
$$
\mbox{LiCoO}_2 + 6\mbox{C} \longleftrightarrow \mbox{Li}_{1-y}\mbox{CoO}_2 +\mbox{Li}_y\mbox{C}_6 \ \ \ \ (0\leq y \leq 1),
$$
$$
\mbox{LiFePO}_4 + 6\mbox{C} \longleftrightarrow \mbox{Li}_{1-y}\mbox{FePO}_4 +\mbox{Li}_y\mbox{C}_6 \ \ \ \ (0\leq y \leq 1),
$$
$$
\mbox{LiMn}_2\mbox{O}_4 + 6\mbox{C} \longleftrightarrow \mbox{Li}_{1-y}\mbox{Mn}_2\mbox{O}_4+\mbox{Li}_y\mbox{C}_6 \ \ \ \ (0\leq y \leq 1).
$$

A 1D electrochemical model is considered for the evolution of the Li concentration $c_{\rm e}(x,t)$ (mol / m$^3$) and the electric potential $\phi_{\rm e}(x,t)$ (V) in the electrolyte and the electric potential
$\phi_{\rm s}(x,t)$ (V) in the solid electrodes, along the $x$--direction, with $x\in (0,L)$ and $L=L_1+\delta +L_2$ being the cell width (m). We assume that $(0,L_{1})$ corresponds to the negative electrode,
($L_1,L_1+\delta)$ corresponds
to the separator and $(L_1+\delta,L)$ corresponds to the positive electrode. This is coupled with a 1D microscopic solid diffusion model for the evolution of the Li concentration $c_{\rm s}(x;r,t)$ in a
generic solid spherical electrode particle (situated at point $x\in (0,L_1)\cup (L_1+\delta , L)$) along the radial $r$--direction, with $r\in [0,R_{\rm s}]$ and $R_{\rm s}$ (m) the average radius of a generic solid active
material particle. This 1D approximation is valid since the characteristic length scale of a typical Li-ion cell along the $x$-axis is on the order of 100 $\mu$m, whereas the characteristic length scale for the remaining two axes
is on the order of 100,000 $\mu$m or more (see \cite{2010ChEtAl}). $R_{\rm s}$ can be different at each electrode and therefore we consider

$$R_{\rm s} = R_{\rm s}(x)=
\left\{
\begin{array}{ll}
 R_{\rm s,-} & \mbox{ if } x\in (0, L_1), \\
 R_{\rm s,+} & \mbox{ if } x\in (L_1+\delta ,L).
\end{array}
\right.
$$

The percentage of available local energy at time $t$ (s) and radius $r$ of a generic negative electrode particle situated at point $x\in (0,L_1)$ in the cell $x$--direction is the the same as its
stoichiometry value $y=y(x,r,t)$, which can be computed as
$$
y(x,r,t)= \frac{c_{\rm s}(x;r,t)}{c_{\rm s,-,max}},
$$
where $c_{\rm s,-,max}$ (mol m$^{-3}$) is the maximum possible concentration in the solid negative electrode and
$c_{\rm s}(x;r,t)$ (mol m$^{-3}$) is the the Li concentration at time $t$, radius $r$ and point $x$. Therefore, the bulk state of charge (SOC) for the negative electrode (it can be also
done for the positive electrode but both are related and, therefore, it suffices to use only one of them) is
\begin{equation} \label{soct}
 \mbox{SOC}(t)=\frac{3}{L_1(R_{\rm s})^3} \int_0^{L_1}\int_0^{R_{\rm s,-}} r^2 \frac{c_{\rm s}(x;r,t)}{c_{\rm s,-,max}}{\rm d}r{\rm d}x.
\end{equation}

We point out that SOC is an nondimensional quantity that can be used as an indicator of the available energy in the cell. Measuring $c_{\rm s}(x;r,t)$ is not easy, but
we can use the mathematical model below (see, e.g., \cite{2010ChEtAl,1993DoFuNe,2006SmWa,2006-b-SmWa}) to compute it. From a theoretical point of view SOC could be a value between 0 and 1 but,
as mentioned above, in practical operating cases it never attains the extreme values 0 or 1. Therefore intermediate values $y_{0\%}$, $z_{0\%}$, $y_{100\%}$, $z_{100\%}$ are used for each
electrode to refer to 0 \% or 100 \% SOC.

\subsection{A full mathematical model} \label{ssafmm}
Based on the models appearing in the literature (see~\cite{1993DoFuNe,1994FuDoNe,2002GoEtAl,2005FaPl,2006SmWa,2006-b-SmWa,2007SmRaWa,2010ChEtAl,2011GuSiWh,2012KiSmIrPe}) and assuming constant diffusion and activity electrolyte coefficients,
a full mathematical model for the performance of a battery, also including heat transfer
dynamics, is given by system of equations (\ref{eqce})--(\ref{eqt}):

\begin{equation} \label{eqce}
\left\{
\begin{array}{l}
{\displaystyle\varepsilon_{\rm e} \frac{\partial  c_{\rm e}}{\partial t} -
D_{\rm e} \frac{\partial }{\partial x} \left( \varepsilon_{\rm e}^p \frac{\partial c_{\rm e}}{\partial x} \right) = \frac{1-t^0_+}{F}
j^{\rm Li}, }
\hspace*{.2cm}\mbox{in } (0,L)\times (0,t_{\rm end}),
\\[.3cm]
{\displaystyle
\frac{\partial c_{\rm e}}{\partial x} (0,t) = \frac{\partial c_{\rm e}}{\partial x} (L,t) =0, \ \ \ t\in (0,t_{\rm end}),
} \\[.3cm]
{\displaystyle
c_{\rm e} (x,0) = c_{\rm e,0} (x), \ \ \ x\in (0,L),
}
\end{array}
\right.
\end{equation}

\begin{equation} \label{eqcs}
 \left\{
\begin{array}{l}
\mbox{For each }x\in (0,L_1) \cup (L_1+\delta,L): \\[.3cm]
{\displaystyle
\frac{\partial  c_{\rm s}}{\partial t} -
\frac{D_{\rm s}}{r^2} \frac{\partial }{\partial r} \left(  r^2 \frac{\partial c_{\rm s}}{\partial r} \right)
= 0, \ \mbox{ in } (0,R_{\rm s})\times (0,t_{\rm end}),
} \\[.3cm]
{\displaystyle \frac{\partial c_{\rm s}}{\partial r} (x;0,t) = 0, \ \  -D_{\rm s}\frac{\partial c_{\rm s}}{\partial r} (x;R_{\rm s},t) =
\frac{R_{\rm s}(x)}{3\varepsilon_{\rm s}(x)F}
j^{\rm Li} , \ \ t\in (0,t_{\rm end}),} \\[.3cm]
 c_{\rm s}(x;r,0)= c_{\rm s,0}(x;r),
\end{array}
\right.
\end{equation}

\begin{equation} \label{eqpe}
 \left\{
\begin{array}{l}
{\displaystyle \mbox{For each }t\in (0,t_{\rm end}):} \\[.3cm]
{\displaystyle
-\frac{\partial }{\partial x} \left( \varepsilon_{\rm e}^p\kappa \frac{\partial \phi_{\rm e}}{\partial x} \right)
+ (1-2t^0_+)\frac{R  T}{F}\frac{\partial }{\partial x} \left( \varepsilon_{\rm e}^p\kappa \frac{\partial}{\partial x}  \ln \big( c_{\rm e}\big) \right)
} = j^{\rm Li} \mbox{ in } (0,L),  \\[.3cm]
{\displaystyle \frac{\partial \phi_{\rm e}}{\partial x} (0,t) = \frac{\partial \phi_{\rm e}}{\partial x} (L,t) = 0},
\end{array}
\right.
\end{equation}

\begin{equation} \label{eqps}
 \left\{
\begin{array}{l}
{\displaystyle \mbox{For each }t\in (0,t_{\rm end}):} \\[.3cm]
{\displaystyle
-\varepsilon_{\rm s} \sigma \frac{\partial^2 \phi_{\rm s}}{\partial x^2}
= } -j^{\rm Li}  \ \mbox{ in } (0,L_1)\cup (L_1+\delta ,L),
 \\[.5cm]
{\displaystyle \varepsilon_{\rm s}(0)\sigma (0) \frac{\partial \phi_{\rm s}}{\partial x} (0,t) = \varepsilon_{\rm s}(L)\sigma (L) \frac{\partial \phi_{\rm s}}{\partial x}
(L,t) = -\frac{I(t)}{A},} \\[.3cm]
{\displaystyle \frac{\partial \phi_{\rm s}}{\partial x} (L_1,t) = \frac{\partial \phi_{\rm s}}{\partial x} (L_1+\delta ,t) = 0}, \\[.3cm]
\end{array}
\right.
\end{equation}

\begin{equation} \label{eqt}
\left\{
\begin{array}{l}
{\displaystyle MC_p \frac{{\rm d} T}{{\rm d}t}} = - h A_{\rm s} \Big( T-T_{\rm amb}\Big) + q_{\rm r}+q_{\rm j}+ q_{\rm c} + q_{\rm e},
 \  t\in (0,t_{\rm end}), \\[.2cm]
T(0)=T_0,
\end{array}
\right.
\end{equation}

In the above system of equations $c_{\rm e}= c_{\rm e}(x,t)$ (mol / m$^3$) at time $t$ (s), $c_{\rm s} = c_{\rm s} (x;r,t)$ (mol / m$^3$), $\phi_{\rm e} = \phi_{\rm e} (x,t)$ (V), $\phi_{\rm s} = \phi_{\rm s} (x,t)$ (V),
$T=T(t)$ is the temperature of the cell (K), $I=I(t)$ is the applied current (A),
$$\varepsilon_{\rm e}=\varepsilon_{\rm e}(x)=
\left\{
\begin{array}{ll}
 \varepsilon_{{\rm e},-} & \mbox{ if } x\in (0, L_1), \\
 \varepsilon_{\rm e,sep} & \mbox{ if } x\in (L_1, L_1+\delta), \\
 \varepsilon_{{\rm e},+} & \mbox{ if } x\in (L_1+\delta ,L)
\end{array}
\right.
$$
is the volume fraction of the electrolyte, $p$ is the Bruggeman porosity exponent (nondimensional constant), $D_{\rm e}$ is the electrolyte diffusion coefficient (m$^2$ s$^{-1}$), $t^0_+$ (dimensionless and assumed here to be constant) is
the transference number of Li$^+$,
$$D_{\rm s} = D_{\rm s} (x) =
\left\{
\begin{array}{ll}
 D_{\rm s,-} & \mbox{ if } x\in (0, L_1), \\
 D_{\rm s,+} & \mbox{ if } x\in (L_1+\delta ,L)
\end{array}
\right.
$$
is the solid phase Li diffusion coefficient (m$^2$ s$^{-1}$),
$$\varepsilon_{\rm s} =\varepsilon_{\rm s}(x)=
\left\{
\begin{array}{ll}
 \varepsilon_{{\rm s},-} & \mbox{ if } x\in (0, L_1), \\
 \varepsilon_{{\rm s},+} & \mbox{ if } x\in (L_1+\delta ,L)
\end{array}
\right.
$$
is the volume fraction of the active materials in the electrodes, $\kappa= \kappa \left(c_{\rm e}(x,t),T(t)\right)$ is the electrolyte phase ionic conductivity (S m$^{-1}$),
$A$ (m$^2$) is the cross-sectional area (also the electrode plate area),
$$\sigma_{\rm s} =\sigma_{\rm s}(x)=
\left\{
\begin{array}{ll}
 \sigma_{{\rm s},-} & \mbox{ if } x\in (0, L_1), \\
 \sigma_{{\rm s},+} & \mbox{ if } x\in (L_1+\delta ,L)
\end{array}
\right.
$$
is the electrical conductivity of solid active materials in an electrode (S m$^{-1}$) and
$$
j^{\rm Li}=j^{\rm Li} \Big( x,\phi_{\rm s}(x,t),\phi_{\rm e}(x,t), c_{\rm s}(x;R_{\rm s}(x),t),c_{\rm e}(x,t),T(t)\Big)
$$
is the reaction current resulting in
 production
 or consumption of Li (A m$^{-3}$) at point $x$ and time $t$. For $j^{\rm Li}$ the Butler-Volmer equation is usually used (see, e.g., \cite{1993DoFuNe,2002GoEtAl,2006SmWa,2006-b-SmWa,2007SmRaWa,2010ChEtAl,2012KiSmIrPe}):
\begin{equation} \label{fcbv}
\hspace*{-.2cm}j^{\rm Li} \Big( x,\phi_{\rm s},\phi_{\rm e}, c_{\rm s}, c_{\rm e},T\Big) =
\left\{
\begin{array}{l}
 {\displaystyle \frac{3\varepsilon_{\rm s} (x) }{R_{\rm s}(x)} i_0  \left[ \exp \left( \frac{\alpha_{\rm a}F}{R \ T}\eta \right) - \exp \left( \frac{-\alpha_{\rm c}F}{R \ T}\eta \right) \right] } \\[.4cm]
 \hspace*{2cm} \mbox{ if } x\in
 (0,L_1)\cup (L_1+\delta , L), \\[.4cm]
 0 \mbox{ if } x\in (L_1,L_1+\delta)
\end{array}
\right.
\end{equation}
(here, for the sake of simplicity, we have considered the solid/electrolyte interfacial film resistance to be zero and therefore is not included in the above equation), where
$a_{\rm s}(x)=\frac{3\varepsilon_{\rm s} (x) }{R_{\rm s}(x)}$ is the specific interfacial area of electrodes (m$^{-1}$), $i_0=i_0 (x,c_{\rm s},c_{\rm e})$,
\begin{equation} \label{ddio}
 i_0 (x,c_{\rm s},c_{\rm e})=
\left\{
\begin{array}{ll}
  k_-(c_{\rm e})^{\alpha_{\rm a}} (c_{\rm s,-,max} -c_{\rm s})^{\alpha_{\rm a}} (c_{\rm s})^{\alpha_{\rm c}} & \mbox{ if } x\in (0, L_1), \\
  k_+(c_{\rm e})^{\alpha_{\rm a}} (c_{\rm s,+,max} -c_{\rm s})^{\alpha_{\rm a}} (c_{\rm s})^{\alpha_{\rm c}} & \mbox{ if } x\in (L_1+\delta ,L)
\end{array}
\right.
\end{equation}
is the exchange current density of an electrode reaction
(A m$^{-2}$), $k_-,k_+$ are kinetic rate constants (A m$^{-2+6\alpha_{\rm a}+3\alpha_{\rm c}}$ mol$^{-2\alpha_{\rm a}-\alpha_{\rm c}}$), $\alpha_{\rm a}$, $\alpha_{\rm c}$ (dimensionless constants) are
anodic and cathodic transfer coefficients for an electrode reaction,
$$
\eta=\eta \Big( x,\phi_{\rm s}(x,t),\phi_{\rm e}(x,t), c_{\rm s}(x;R_{\rm s}(x),t),T(t)\Big),
$$
$$
\eta\Big( x,\phi_{\rm s},\phi_{\rm e}, c_{\rm s},T\Big) = \left\{
\begin{array}{l}
 {\displaystyle \phi_{\rm s} - \phi_{\rm e} - U(x,c_{\rm s},T), \mbox{ if } x\in
 (0,L_1)\cup (L_1+\delta , L),} \\[.2cm]
 0 \mbox{ if } x\in (L_1,L_1+\delta)
\end{array}
\right.
$$
is the surface overpotential (V) of an electrode reaction and
$U(x,c_{\rm s},T)$ is the equilibrium potential (V) at the solid/electrolyte interface (i.e. the open circuit voltage - OCV). A way of expressing $U$ is (see \cite{2011GuSiWh}):

$$ 
 U(x,c,T)=
$$
$$
\left\{
\begin{array}{ll}
 {\displaystyle U_{-}\left( \frac{c}{c_{\rm s,-,max}}\right) + \frac{\partial U_-}{\partial T}\left( \frac{c}{c_{\rm s,-,max}}\right)\Big( T-T_{\rm ref}\Big) }
 & \mbox{ if } x\in (0, L_1), \\[.4cm]
 {\displaystyle U_{+}\left( \frac{c}{c_{\rm s,+,max}}\right) + \frac{\partial U_+}{\partial T}\left( \frac{c}{c_{\rm s,+,max}}\right)\Big( T-T_{\rm ref}\Big) }
 & \mbox{ if } x\in (L_1+\delta ,L),
\end{array}
\right.
$$
where ${\displaystyle U_{-},U_{+}, \frac{\partial U_-}{\partial T} \frac{\partial U_-}{\partial T}}$ are functions tipically obtained from fitting experimental data
and $c_{\rm s,+,max}$ (mol m$^{-3}$) is the maximum possible concentration in
the solid positive electrode.

\begin{remark} \label{reliap}
 For low overpotential cases Butler-Volmer equation (\ref{fcbv}) can be simplified to the following linearized version (see \cite{2007SmRaWa}):

$$
\hspace*{-.2cm}j^{\rm Li} \Big( x,\phi_{\rm s},\phi_{\rm e}, c_{\rm s}, c_{\rm e},T\Big) =
\left\{
\begin{array}{l}
 {\displaystyle \frac{3\varepsilon_{\rm s} (x) }{R_{\rm s}(x)} i_0 \frac{(\alpha_{\rm a}+\alpha_{\rm c})F}{R \ T}\eta },  \mbox{ if } x\in
 (0,L_1)\cup (L_1+\delta , L), \\[.4cm]
 0, \mbox{ if } x\in (L_1,L_1+\delta).
\end{array}
\right.
$$
\end{remark}

\vspace{.3cm}

Regarding heat transfer Eq.~(\ref{eqt}), $M$ (kg) is the mass of the battery, $C_p$ (J kg$^{-1}$ K$^{-1}$) is the specific heat capacity, $h$ (W m$^{-2}$ K$^{-1}$) is the heat transfer
coefficient for convection, $A_{\rm s}$ (m$^2$) is the cell surface
area exposed to the convective cooling medium (typically air), $T_{\rm amb}$ is the (ambient) temperature of the cooling medium,
$$
q_{\rm r}= q_{\rm r} (t) = A\int_0^L j^{\rm Li} \eta  {\rm d}x
$$
is the total reaction heat, $q_{\rm j} =q_{\rm j} (t)$,
\begin{eqnarray*}
q_{\rm j} (t) & = & {\displaystyle A\varepsilon_{\rm s,-}\sigma_-\int_0^{L_1}  \left( \frac{\partial \phi_{\rm s}}{\partial x} \right)^2 {\rm d}x +
A\varepsilon_{\rm s,+}\sigma_+ \int_{L_1+\delta}^L  \left( \frac{\partial \phi_{\rm s}}{\partial x} \right)^2 {\rm d}x} \\
& & {\displaystyle +   A\int_0^L \left[  \varepsilon_{\rm e}^p\kappa  \left( \frac{\partial \phi_{\rm e}}{\partial x} \right)^2   + (2t^0_+ - 1)\frac{R \ T}{F} \varepsilon_{\rm e}^p\kappa
\left( \frac{\partial \ln c_{\rm e}}{\partial x}  \right)
\left( \frac{\partial \phi_{\rm e}}{\partial x} \right) \right] {\rm d}x }
\end{eqnarray*}
is the ohmic heat due to the current carried in each phase and the limited conductivity of that phase,
$$
q_{\rm c}= q_{\rm c}(t)= I(t)^2\frac{R_{\rm f}}{A}
$$
is the ohmic heat generated in the cell due to contact resistance between current collectors and electrodes, $R_{\rm f}$ ($\Omega$ m$^2$) is the film resistance of the electrodes and
$$
\begin{array}{lll}
{\displaystyle  q_{\rm e} = q_{\rm e}(t)}  & = & {\displaystyle TA\left[  \int_0^{L_1} j^{\rm Li}
\frac{\partial U_{-}}{\partial T} \left( \frac{c_{\rm s}(x;R_{\rm s,-},t)}{c_{\rm s,-,max}}\right) {\rm d}x  \right. }   \\[.4cm]
& & {\displaystyle \left. + \int_{L_1+\delta}^L j^{\rm Li}
\frac{\partial U_{+}}{\partial T} \left( \frac{c_{\rm s}(x;R_{\rm s,+},t)}{c_{\rm s,+,max}}\right) {\rm d}x \right]  }
\end{array}
$$
is the reversible heat caused by the reaction entropy change (see \cite{2001ThEtAl} and \cite{2011GuSiWh} for a simpler formulation of $q_{\rm e}(t)$).

Heat sources $q_{\rm r}$, $q_{\rm j}$ and $q_{\rm c}$ are always positive and, as explained in \cite{2006SmWa}, the second term inside the last integral of $q_{\rm j}$ is generally
negative. On the other hand, $q_{\rm e}$ can be either positive or negative.

\begin{remark} The second term on the left hand side of system (\ref{eqpe}) is often written in the literature using $(1-t^0_+)$ (see Refs.~\cite{1993DoFuNe} and \cite{1994FuDoNe}) or $2(1-t^0_+)$
(see Refs.~\cite{2002GoEtAl,2006SmWa,2006-b-SmWa,2007SmRaWa,2010ChEtAl,2012KiSmIrPe}), instead of $(1-2t^0_+)$, which is the correct term (see pag. 250 of Ref.~\cite{1973Ne}).
We remind that we are assuming here that $t^0_+$ is constant and it is given (see \cite{1973Ne}) by
$$
t^0_+= \frac{D_{0,+}}{D_{0,-}+D_{0,+}},
$$
with $D_{0,-}, D_{0,+}$ transport coefficients related to the diffusion of Li$^+$ cations and the corresponding anions, respectively.
\end{remark}

\begin{remark} Many authors (see, e.g., \cite{2006SmWa} and \cite{2011GuSiWh}) use some physiochemical parameters involved in the model depending on $T$ and use
the Arrhenius equation to express such a dependency, which needs the corresponding activation energies. The dependency of some parameters on $T$ can sometimes be estimated. For instance
in Ref.~\cite{2012KiSmIrPe} the authors give formulae for $\kappa(c_{\rm e},T), D_{\rm e}(c,T)$ and $t_+^0(c,T)$.
\end{remark}

\begin{remark} \label{rsfp} This system of equations does not have uniqueness of solution (if $\phi_{\rm s}(x,t)$ and $\phi_{\rm e}(x,t)$ are solutions
then $\phi_{\rm s}(x,t)+c(t)$ and $\phi_{\rm e}(x,t)+c(t)$ are also solutions, for any function $c(t)$). A way of avoiding that is to set a reference value of
$\phi_{\rm s}(x,t)$ or $\phi_{\rm e}(x,t)$ at some point $x$.
For instance we can impose $\phi_{\rm s} (0,t)=0$ for any $t\in [0,t_{\rm end}]$ (see Remark \ref{eaude}). Some results regarding the existence and uniqueness of solution are showed in Sec.~\ref{eusol}.
\end{remark}

\vspace{.3cm}

After solving the above model, we can get SOC$(t)$ by using (\ref{soct}) and the cell voltage $V(t)$ (A), at time $t$, given by
\begin{equation} \label{erci}
 V(t) = \phi_{\rm s} (L,t) - \phi_{\rm s} (0,t) - \frac{R_{\rm f}}{A}I(t).
\end{equation}
The battery pack voltage is calculated by multiplying the single cell voltage of Eq.~(\ref{erci}) by the number of serially connected cells in the battery.

\section{Correct boundary conditions in (\ref{eqps})} \label{cbcs}

\begin{remark} \label{rmswbc}
In \cite{2006SmWa,2006-b-SmWa,2007SmRaWa,2012KiSmIrPe},
authors use the following boundary conditions at $x=0$ and $x=L$, instead of those presented in (\ref{eqps}):
\begin{equation} \label{wbc1}
 -\varepsilon_{\rm s,-}\sigma_- \frac{\partial \phi_{\rm s}}{\partial x} (0,t) = \varepsilon_{\rm s,+}\sigma_+ \frac{\partial \phi_{\rm s}}{\partial x} (L,t) = \frac{I(t)}{A}.
\end{equation}

In \cite{2002GoEtAl} authors use the following boundary conditions at $x=0$ and $x=L$:
\begin{equation} \label{wbc2}
 -\varepsilon_{\rm s,-}\sigma_- \frac{\partial \phi_{\rm s}}{\partial x} (0,t) = \varepsilon_{\rm s,+}\sigma_+ \frac{\partial \phi_{\rm s}}{\partial x} (L,t) = -\frac{I(t)}{A}.
\end{equation}

In Sec.~\ref{pdedu} it is proved that both conditions are not correct from a physical point of view, because the resulting system of equations
does not model a battery performance. Then, in Sec.~\ref{prelim} it is proved that they are also mathematically incorrect, because, in fact, the corresponding system of boundary value
problems modelling the conservation of charge does not have any
solution (therefore the system is not well-posed) unless $I(t) \equiv 0$ (see Remark \ref{nosolu}).
\end{remark}

\subsection{Deduction of system (\ref{eqps})} \label{pdedu}

Given an external current $I(t)$ (A) applied to the battery ($I(t)>0$ when the battery is discharging), according to Kirchoff's law $I=I_{\rm s}+I_{\rm e}$, where $I_{\rm s},I_{\rm e}$ are the current (A) in the
solid electrode and in the electrolyte, respectively.
Now, by Ohm's law,
\begin{equation} \label{edli}
I_{\rm s}(x,t)=-A\sigma (x) \frac{\partial \phi_{\rm s}}{\partial x}(x,t)=I(t)-I_{\rm e}(x,t).
\end{equation}

\begin{figure}[htb]
\begin{center}
\begin{minipage}[htb]{11cm}
\includegraphics[scale=.6]{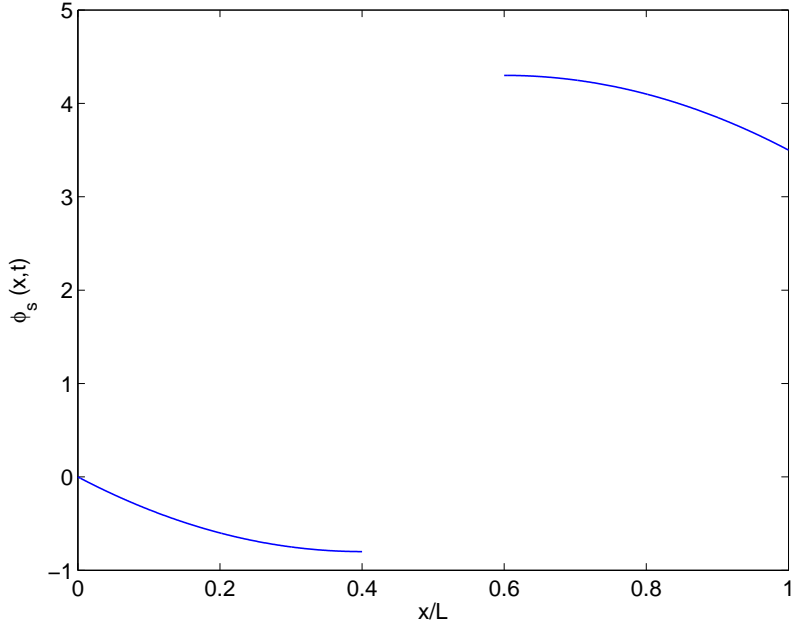}
\caption{\protect\small Typical Graph of $\phi_{\rm s}(\cdot,t)$
for some time $t\in (t,t_{\rm end})$.}
\label{GraphTypPhis}
\end{minipage}
\end{center}
\end{figure}

Then, for each $t\in (0,t_{\rm end})$, taking into account the porosity nature of the solid electrodes (only a fraction of its volume contribute to its conductivity) and using that
\begin{equation} \label{reij}
\frac{\partial I_{\rm e}}{\partial x} (x,t)= A j^{\rm Li},
\end{equation}
the following equation for the conservation of charge in the electrode solid phase is obtained:

$$
-\varepsilon_{\rm s}\sigma \frac{\partial^2 \phi_{\rm s}}{\partial x^2}
= - j^{\rm Li} , \ \mbox{ in } (0,L_1)\cup (L_1+\delta ,L).
$$
We point out that $I_{\rm s}(x,t)=0$ (and therefore $I_{\rm e}(x,t)=I(t)$) for all $x\in (L_1,L_1+\delta)$ and $I_{\rm e}(0,t)=I_{\rm e}(L,t)=0$. This provides, using (\ref{edli}),
the following boundary conditions:
$$
\varepsilon_{\rm s}(0)\sigma (0) \frac{\partial \phi_{\rm s}}{\partial x} (0,t) = \varepsilon_{\rm s}(L)\sigma (L) \frac{\partial \phi_{\rm s}}{\partial x} (L,t) = -\frac{I(t)}{A},
$$
$$
\frac{\partial \phi_{\rm s}}{\partial x} (L_1,t) = \frac{\partial \phi_{\rm s}}{\partial x} (L_1+\delta ,t) = 0,
$$
which completes the deduction of system (\ref{eqps}). Therefore, if for instance the battery is discharging, $j^{\rm Li}>0$ (respectively, $j^{\rm Li}<0$) if $x\in (0,L_1)$ (respectively, if $x\in (L_1+\delta, L)$) and
the graph of $\phi_{\rm s}(\cdot,t)$
for some time $t\in (t,t_{\rm end})$ could be qualitatively (not necessarily quantitatively) similar to that showed in Figure \ref{GraphTypPhis}.

\subsection{Deduction of system (\ref{eqpe})} \label{pdedu2}

Neglecting deviations of the electrolyte solution from ideal behavior
$$
I_{\rm e}(x,t)= -A\kappa \frac{\partial \phi_{\rm e}}{\partial x}(x,t) + A
(1-2t^0_+)\frac{RT}{F} \kappa  \frac{\partial}{\partial x}  \ln \big( c_{\rm e}(x,t)\big)
$$
(see \cite{1973Ne}). The first term on the right hand side is due to Ohm's law (as in (\ref{edli})) and the second term accounts for concentration variations.

Then, for each $t\in (0,t_{\rm end})$, taking into account the porosity nature of the solid electrodes and using again (\ref{reij}),
the following equation for the conservation of charge in the electrode solid phase is obtained:
$$
-\frac{\partial }{\partial x} \left( \varepsilon_{\rm e}^p\kappa  \frac{\partial \phi_{\rm e}}{\partial x} \right)
+ (1-2t^0_+)\frac{RT}{F}\frac{\partial }{\partial x} \left( \varepsilon_{\rm e}^p\kappa  \frac{\partial}{\partial x}  \ln \big( c_{\rm e}\big) \right)
= j^{\rm Li}  , \ \mbox{ in } (0,L),
$$
together with boundary conditions $I_{\rm e}(0,t)=I_{\rm e}(L,t)=0$. i.e.,
$$
-\kappa (c_{\rm e}(0,t)) \frac{\partial \phi_{\rm e}}{\partial x} (0,t)
+ (1-2t^0_+)\frac{2T}{F} \kappa (c_{\rm e}(0,t))  \frac{\partial}{\partial x}  \ln \big( c_{\rm e}(0,t)\big) =0,
$$
$$
-\kappa (c_{\rm e}(L,t)) \frac{\partial \phi_{\rm e}}{\partial x} (L,t)
+ (1-2t^0_+)\frac{2T}{F} \kappa (c_{\rm e}(L,t))  \frac{\partial}{\partial x}  \ln \big( c_{\rm e}(L,t)\big) =0,
$$
which, taking into account the boundary condtions in (\ref{eqce}), are equivalent to
$$
\frac{\partial \phi_{\rm e}}{\partial x} (0,t) = \frac{\partial \phi_{\rm e}}{\partial x} (L,t) = 0.
$$

\section{About existence and uniqueness of solution of system (\ref{eqpe})--(\ref{eqps})} \label{eusol}

Let us study some properties regarding the existence and uniqueness of solution of (\ref{eqpe})--(\ref{eqps}).
We will only focus here on the equations involving conservation of charge.

Even if coupling with the concentrations of Lithium is not included in the analytic work, the system under studied is interesting and not trivial at all. It is interesting because suitable numerical schemes of the full model will probably need to solve this kind of systems at each time step of a time discretized version of the full model. It is non trivial because it is a system of three boundary value problems, defined and coupled over 3 different domains. Moreover, as we will see below, the system does not have uniqueness of solution except if the parameters involved in the system satisfy a compatibility condition (which is satisfied when the parameters come from the full model) and we add an extra boundary condition (which is equivalent to set a reference potencial value and, therefore, consistent with the physics of the problem). The existence result showed in Section 4.3 is not standard, since it requires the use of Functional Analysis techniques.

\subsection{Preliminaries} \label{prelim}

Let us suppose $(c_{\rm e}, c_{\rm s}, \phi_{\rm e}, \phi_{\rm s}, T)$ is a solution of the full model. Given
$t\in \left( 0,t_{\rm end}\right)$ let us denote $u(x)=\phi_{\rm e}(x,t)$ and $v(x)=\phi_{\rm s}(x,t)$. Then, if
$$
E(x,w)=\exp\left( \frac{\alpha_{\rm a}F}{RT} \big( w-f(x)\big) \right) - \exp\left( \frac{-\alpha_{\rm c}F}{RT} \big( w-f(x)\big) \right) ,
$$
$u$ and $v$ must satisfy the following system of equations equivalent to
(\ref{eqpe})--(\ref{eqps}):

\begin{equation} \label{ndepef}
 \left\{
\begin{array}{l}
{\displaystyle
- \frac{{\rm d} }{{\rm d} x} \left( r (x) \frac{{\rm d} u}{{\rm d} x} (x) \right)  + \nu \frac{{\rm d} }{{\rm d} x} \left( r(x) \frac{{\rm d} }{{\rm d} x} \ln (c(x)) \right)
} \\
\hspace*{2cm}= \left\{
\begin{array}{l}
 {\displaystyle k_-b_-(x) E(x,v(x)-u(x)) , \mbox{ if } x\in (0,L_1),} \\[.3cm]
 {\displaystyle k_+b_+(x) E(x,v(x)-u(x)) , \mbox{ if }  x\in (L_1+\delta ,L),} \\[.3cm]
 0 , \mbox{ if }  x\in (L_1,L_1+\delta),
\end{array}
\right.  \\[.5cm]
{\displaystyle \frac{{\rm d} u}{{\rm d} x} (0) = \frac{{\rm d} u}{{\rm d} x} (L) = 0},
\end{array}
\right.
\end{equation}

\begin{equation} \label{ndepsf}
 \left\{
\begin{array}{l}
{\displaystyle
-\frac{{\rm d}^2 v}{{\rm d} x^2}(x)
= } {\displaystyle -q_-b_-(x)E(x,v(x)-u(x)) , \mbox{ if }   x\in (0,L_1),}
 \\[.5cm]
 {\displaystyle
-\frac{{\rm d}^2 v}{{\rm d} x^2} (x)
= } {\displaystyle -q_+b_+(x)E(x,v(x)-u(x)) , \mbox{ if } x\in (L_1+\delta ,L) ,}
 \\[.5cm]
{\displaystyle - \frac{{\rm d} v}{\partial x} (0) = d_{-}, \ \ \
\frac{{\rm d} v}{\partial x}
(L) =- d_+ , \ \
 \frac{{\rm d} v}{{\rm d} x} \left(L_1\right) =
 \frac{{\rm d} v}{{\rm d} x} \left(L_1+\delta \right) = 0},
\end{array}
\right.
\end{equation}
$$r(x)= \varepsilon_{\rm e}(x)^p \kappa (c_{\rm e}(x,t),T(t)), \ c(x)=c_{\rm e}(x,t), \ \nu = (1-2t^0_+)\frac{2 T(t)}{F},
$$
$$
b_\pm(x)=
  \frac{3\varepsilon_{\rm s} (x) }{R_{\rm s}(x)}(c_{\rm e}(x,t))^{\alpha_{\rm a}} (c_{\rm s,\pm,max} -c_{\rm s}(x;R_{\rm s},t))^{\alpha_{\rm a}} (c_{\rm s}(x;R_{\rm s},t))^{\alpha_{\rm c}},
$$
\begin{equation} \label{ddpa}
 q_\pm = \frac{k_\pm}{\varepsilon_\pm \sigma_\pm} , \
d_{\pm}=\frac{1}{\varepsilon_{\rm s,\pm}\sigma_\pm }\frac{I(t)}{A}, \ f(x)= U(x,c_{\rm s}(x;R_{\rm s},t),T(t)).
\end{equation}

\begin{remark} Applying the divergence theorem it is easy to show that any solution of (\ref{ndepef})--(\ref{ndepsf}) must satisfy
$$
k_-\int_0^{L_1} b_-(x) E(x,v(x)-u(x)) {\rm d}x
+k_+\int_{L_1+\delta}^1 b_+(x) E(x,v(x)-u(x)) {\rm d}x =0,
$$
$$
q_-\int_0^{L_1} b_-(x) E(x,v(x)-u(x)) {\rm d}x = d_- ,
$$
$$
q_+\int_{L_1+\delta}^{L} b_+(x) E(x,v(x)-u(x)) {\rm d}x = -d_+.
$$
Therefore, a necessary compatibility condition for the existence of solution of system (\ref{ndepef})-(\ref{ndepsf}) is $\frac{k_-d_-}{q_-}=\frac{k_+d_+}{q_+}$,
which is true if (\ref{ddpa}) holds, because $\frac{k_\pm d_\pm}{q_\pm}= \frac{I(t)}{A}$.
\end{remark}

\begin{remark} \label{nosolu}
 If we use boundary conditions (\ref{wbc1}) or (\ref{wbc2}) the corresponding compatibility condition would be  $\frac{k_-d_-}{q_-}=-\frac{k_+d_+}{q_+}$, or equivalently $I(t)=-I(t)$, which is only
satisfied if $I(t)\equiv 0$.
\end{remark}

Since function $r(x)$ may not be smooth (at least in a general situation), system (\ref{ndepef})--(\ref{ndepsf}) (or, equivalently, (\ref{eqpe})--(\ref{eqps})) may not have a classical
 solution.
 Therefore, we look for what is commonly known as a {\em weak solution}. In order to find the requirements that this solution shoud satisfy, we use the following functional
 analysis framework.

Let us consider the Hilbert space
$$H=H^1(0,L)\times \Big( H^1(0,L_1)\cap H^1(L_1+\delta ,L)\Big),
$$
with the norm given by
\begin{eqnarray*}
\| (u,v)\|_V^2 & = &  \| u\|^2_{H^1(0,L)} + \| v\|^2_{H^1(0,L_1)} + \| v\|^2_{H^1(L_1+\delta,L)} \\
& = & {\displaystyle \int_0^L \left[ u(x)^2 + \left( \frac{{\rm d} u}{{\rm d} x} (x)\right)^2 \right] {\rm d}x } \\
& & {\displaystyle \hspace*{-1cm} + \int_0^{L_1} \left[ v(x)^2 + \left( \frac{{\rm d} v}{{\rm d} x} (x)\right)^2 \right] {\rm d}x
+ \int_{L_1+\delta}^1 \left[ v(x)^2 + \left( \frac{{\rm d} v}{{\rm d} x} (x)\right)^2 \right] {\rm d}x .
}
\end{eqnarray*}
Here
$$
H^1(a,b)=\left\{ \phi \in L^2(a,b): \frac{{\rm d} \phi}{{\rm d} x}\in L^2(a,b)\right\},
$$
where $L^p(a,b)$ (with $p\geq 1$) is the set of all measurable functions from $(a,b)$ to $\R$ whose absolute value raised to the $p$-th power has finite integral and $ \frac{{\rm d} \phi}{{\rm d} x}$ denotes
the derivative of $\phi$ in the sense of distributions.

Multiplying (\ref{ndepef}) (duality product $<\cdot , \cdot >_{(H^1(0,L))^\prime \times H^1(0,L)}$) by a function $\varphi\in H^1(0,L)$ and similarly in (\ref{ndepsf}) by a function
$\psi\in H^1(0,L_1)\cap H^1(L_1+\delta ,L)$ we get
$$
\int_0^L r(x) \frac{{\rm d} u}{{\rm d} x} (x) \frac{{\rm d} \varphi }{{\rm d} x} (x) {\rm d}x
+ \int_0^{L_1} \frac{{\rm d} v}{{\rm d} x} (x) \frac{{\rm d} \psi }{{\rm d} x} (x) {\rm d}x
+ \int_{L_1+\delta}^{L} \frac{{\rm d} v}{{\rm d} x} (x) \frac{{\rm d} \psi }{{\rm d} x} (x) {\rm d}x
$$
$$
-k_-\int_0^{L_1} b_-(x) E(x,v(x)-u(x)) \varphi (x) {\rm d}x
-k_+\int_{L_1+\delta}^L b_+(x) E(x,v(x)-u(x)) \varphi (x) {\rm d}x
$$
$$
+q_-\int_0^{L_1} b_-(x) E(x,v(x)-u(x)) \psi (x) {\rm d}x
+q_+\int_{L_1+\delta}^L b_+(x) E(x,v(x)-u(x)) \psi (x) {\rm d}x
$$
$$
 =\nu \int_0^L \frac{r(x)}{c(x)} \frac{{\rm d} c}{{\rm d} x} (x) \frac{{\rm d} \varphi }{{\rm d} x} (x) {\rm d}x - d_- \psi (0) +d_+\psi (1).
$$
This equality can be rewritten as
$$
a\left[ (u,v) \ , (\varphi ,\psi) \right] = l \left[ (\varphi, \psi)\right] ,
$$
where $a\left[ (u,v) \ , (\varphi ,\psi) \right]$ and $l\left[ (\varphi, \psi)\right]$ are the left hand side and right hand side, respectively, of the previous expression. According
to this, we have the following definition of a weak solution (that we will call {\em solution} in the following).

\begin{definition}
 A (weak) solution of system (\ref{ndepef})--(\ref{ndepsf}) is a couple of functions $(u,v) \in H$ such that
\begin{equation} \label{edws}
 a\left[ (u,v) \ , (\varphi ,\psi) \right] = l\left[ (\varphi, \psi)\right] , \ \forall \ (\varphi , \psi)\in H.
\end{equation}
\end{definition}

\subsection{About the uniqueness of solution of systems  (\ref{ndepef})--(\ref{ndepsf})}

For simplicity let us study system (\ref{ndepef})--(\ref{ndepsf}), which is equivalent to system (\ref{eqpe})--(\ref{eqps}).

\begin{proposition} \label{tsuec}
 If $(u_1,v_1)$, $(u_2,v_2)$ are two solutions of (\ref{ndepef})--(\ref{ndepsf}), then there exists a constant $s\in \R$ such that
 $$
 u_2(x)-u_1(x) \equiv s, \ \forall \ x\in (0,L).
 $$
and
$$
 v_2(x)-v_1(x) \equiv s, \ \forall \ x\in \left( 0,L_1\right) \cap \left( L_1+\delta ,L\right).
$$
\end{proposition}

\noindent
{\bf Proof:}
 Let us suppose $(u_1,v_1)$, $(u_2,v_2)$ are two solutions of (\ref{ndepef})--(\ref{ndepsf}), i.e., they satisfy (\ref{edws}). Then, using $(\varphi , \psi)= (u_2-u_1,
\overline{v}_2-\overline{v}_1)$ in
(\ref{edws}), with
$$
\overline{v}_i(x)= \left\{
\begin{array}{ll}
 {\displaystyle \frac{k_-}{q_-}v_i(x)} & {\displaystyle \mbox{ if } x\in \left( 0, L_1\right) , } \\[.3cm]
 {\displaystyle \frac{k_+}{q_+}v_i(x)} & {\displaystyle \mbox{ if } x\in \left(L_1+\delta , L\right) ,}
\end{array}
\right. \ (i=1,2),
$$
we have
$$
 a\left[ (u_2,v_2) \ , (u_2-u_1,\overline{v}_2-\overline{v}_1) \right] - a\left[ (u_1,v_1) \ , (u_2-u_1,\overline{v}_2-\overline{v}_1) \right]= 0,
$$
which is equivalent to

$$
\int_0^L r(x) \left( \frac{{\rm d} (u_2-u_1)}{{\rm d} x} (x) \right)^2 {\rm d}x
$$
$$
+ \frac{k_-}{q_-} \int_0^{L_1} \left( \frac{{\rm d} (v_2-v_1)}{{\rm d} x} (x) \right)^2 {\rm d}x
+ \frac{k_+}{q_+} \int_{L_1+\delta}^L \left( \frac{{\rm d} (v_2-v_1)}{{\rm d} x} (x) \right)^2 {\rm d}x
$$
$$
 + k_-\int_0^{L_1} b_-(x) \Big(  E(x,v_2(x)-u_2(x)) \hspace{7cm}
$$
$$\hspace*{1cm} -  E(x,v_1(x)-u_1(x)) \Big)
\Big( \big( v_2(x)-u_2(x)\big) - \big( v_1(x)-u_1(x)\big) \Big) {\rm d}x
$$
$$
 + k_+\int_{L_1+\delta}^L b_+(x) \Big(  E(x,v_2(x)-u_2(x)) \hspace{7cm}
$$
$$\hspace*{1cm} -  E(x,v_1(x)-u_1(x)) \Big)
\Big( \big( v_2(x)-u_2(x)\big) - \big( v_1(x)-u_1(x)\big) \Big) {\rm d}x
 =0.
$$
Hence, from the monotonicity properties of the exponential function,
$$
\int_0^L r(x) \left( \frac{{\rm d} (u_2-u_1)}{{\rm d} x} (x) \right)^2 {\rm d}x
$$
$$
+ \frac{k_-}{q_-} \int_0^{L_1} \left( \frac{{\rm d} (v_2-v_1)}{{\rm d} x} (x) \right)^2 {\rm d}x
+ \frac{k_+}{q_+} \int_{L_1+\delta}^L \left( \frac{{\rm d} (v_2-v_1)}{{\rm d} x} (x) \right)^2 {\rm d}x =0
$$
and
$$
\Big( v_2(x)-u_2(x)\big) - \big( v_1(x)-u_1(x)\Big) =0 \ \ \forall \ x\in \left( 0,L_1\right) \cap \left( L_1+\delta ,L\right) ,
$$
which implies that there exists a constant $s\in \R$ such that
 $$
 u_2(x)-u_1(x) \equiv s, \ \forall \ x\in (0,L)
 $$
and
$$
 v_2(x)-v_1(x) \equiv s, \ \forall \ x\in \left( 0,L_1\right) \cap \left( L_1+\delta ,L\right).
$$

\begin{corollary} \label{csuec}
 If $(u_1,v_1)$, $(u_2,v_2)$ are two solutions of (\ref{ndepef})--(\ref{ndepsf}), then $v_1(x)-u_1(x)=v_2(x)-u_2(x)$ for all $x\in \left( 0,L_1\right) \cap \left( L_1+\delta ,L\right)$ and, therefore,
 $$E(x,v_1(x)-u_1(x))=E(x,v_2(x)-u_2(x)).
 $$
\end{corollary}

\begin{corollary} \label{cospd} If we add the boundary condition
\begin{equation} \label{cfp}
v (0)=0
\end{equation}
(or $v(\overline{x})=s\in\R$ or $u(\overline{x})=s\in\R$, with $\overline{x}$ and $s$ arbitrarily chosen), then if system (\ref{ndepef}), (\ref{ndepsf}), (\ref{cfp}) has a solution, it is unique.
\end{corollary}

\begin{remark} \label{eaude} According to the mathematical results showed in Proposition \ref{tsuec}, Corollary \ref{csuec} and Corollary \ref{cospd}, a way of getting
the property of uniqueness of solution of system (\ref{eqpe})--(\ref{eqps}) (see non-uniqueness results in Remark \ref{rsfp}) is setting a reference potential value for $\phi_{\rm s}$ or $\phi_{\rm e}$ for each $t\in (0,t_{\rm end})$. For instance
\begin{equation} \label{uncon}
 \phi_{\rm s} (0,t)=0  \ \ \mbox{ (or }\phi_{\rm s} (\overline{x},t)=s(t)\in\R, \overline{x} \mbox{ and }s(t) \mbox{ arbitrarily chosen).}
\end{equation}
\end{remark}

\subsection{About the existence of solution of system (\ref{ndepef}), (\ref{ndepsf}), (\ref{cfp}))}

Let us study the existence of solution for the following linearized approximation (see Remark \ref{reliap}) of system (\ref{ndepef}), (\ref{ndepsf}), (\ref{cfp}):

\begin{equation} \label{ndepefl}
 \left\{
\begin{array}{l}
{\displaystyle
- \frac{{\rm d} }{{\rm d} x} \left( r(x) \frac{{\rm d} u}{{\rm d} x} (x) \right)  + \nu \frac{{\rm d} }{{\rm d} x} \left( r(x) \frac{{\rm d} }{{\rm d} x} \ln (c(x)) \right)
} \\
\hspace*{1.5cm}= \left\{
\begin{array}{l}
 {\displaystyle \lambda_- b_-(x) \left( v(x)-u(x)-f(x) \right) , \mbox{ if } x\in (0,L_1),} \\[.3cm]
 {\displaystyle \lambda_+ b_+(x) \left( v(x)-u(x)-f(x) \right) , \mbox{ if } x\in (L_1+\delta ,L),} \\[.3cm]
  0 , \mbox{ if }  x\in (L_1,L_1+\delta),
\end{array}
\right.  \\[.5cm]
{\displaystyle \frac{{\rm d} u}{{\rm d} x} (0) = \frac{{\rm d} u}{{\rm d} x} (L) = 0},
\end{array}
\right.
\end{equation}

\begin{equation} \label{ndepsfl}
 \left\{
\begin{array}{l}
{\displaystyle
-\frac{{\rm d}^2 v}{{\rm d} x^2} (x)
= } {\displaystyle -\theta_- b_-(x) \left( v(x)-u(x)-f(x) \right) , \ x\in (0,L_1),}
 \\[.5cm]
 {\displaystyle
-\frac{{\rm d}^2 v}{{\rm d} x^2} (x)
= } {\displaystyle -\theta_+ b_+(x) \left( v(x)-u(x)-f(x) \right) , \ x\in (L_1+\delta ,L) ,}
 \\[.5cm]
{\displaystyle - \frac{{\rm d} v}{{\rm d} x} (0) = d_{-}, \ \ \
\frac{{\rm d} v}{{\rm d} x}
(1) =- d_+ , \ \ \ \ \ \ \ \
\frac{{\rm d} v}{{\rm d} x} \left(L_1\right) =
\frac{{\rm d} v}{{\rm d} x} \left(L_1+\delta \right) = 0},
\end{array}
\right.
\end{equation}

\begin{equation} \label{cfpl}
 v (0)=0  \ \ \mbox{ (or }v(\overline{x})=s\in\R, \mbox{ with }\overline{x} \mbox{ and }s \mbox{ arbitrarily chosen),}
\end{equation}
where
$$
\lambda_\pm = k_\pm \frac{(\alpha_{\rm a} + \alpha_{\rm c}) F}{RT}, \ \theta_\pm = q_\pm \frac{(\alpha_{\rm a} + \alpha_{\rm c}) F}{RT}.
$$

\begin{remark} \label{remseul}
The results in Proposition \ref{tsuec} are also valid for the solutions of (\ref{ndepefl})--(\ref{ndepsfl}).
Therefore, if (\ref{ndepefl})--(\ref{cfpl}) has a solution, then it is unique.
\end{remark}

Multiplypling in (\ref{ndepefl}) (duality product $<\cdot , \cdot >_{(H^1(0,L))^\prime \times H^1(0,L)}$) by a function $\varphi\in H^1(0,L)$ and similarly in the first and second equations of
(\ref{ndepsfl}) by functions ${\displaystyle \frac{\lambda_-}{\theta_-}\psi}$ and ${\displaystyle \frac{\lambda_+}{\theta_+}\psi}$, respectively, with
$\psi\in H^1(0,L_1)\cap H^1(L_1+\delta ,L)$, the following equality is obtained:
$$
 \int_0^L r(x) \frac{{\rm d} u}{\partial x} (x) \frac{{\rm d} \varphi }{{\rm d} x} (x) {\rm d}x
+ \frac{\lambda_-}{\theta_-}\int_0^{L_1} \frac{{\rm d} v}{{\rm d} x} (x) \frac{{\rm d} \psi }{{\rm d} x} (x) {\rm d}x
+ \frac{\lambda_+}{\theta_+}\int_{L_1+\delta}^L \frac{{\rm d} v}{{\rm d} x} (x) \frac{{\rm d} \psi }{{\rm d} x} (x) {\rm d}x
$$
$$
+\lambda_- \int_0^{L_1} b_-(x)\left( v(x)-u(x)\right) \left( \psi (x) -\varphi (x) \right) {\rm d}x
$$
$$
+\lambda_+ \int_{L_1+\delta}^L b_+(x) \left( v(x)-u(x)\right) \left( \psi (x) -\varphi (x) \right) {\rm d}x
$$

$$
 =  \lambda_- \int_0^{L_1} b_-(x) f(x) \left( \varphi (x) -\psi(x)\right) {\rm d}x
 + \lambda_+ \int_{L_1+\delta}^L b_+(x) f(x) \left( \varphi (x) - \psi(x)\right) {\rm d}x
$$
$$
 + \nu \int_0^L \frac{r(x)}{c(x)} \frac{{\rm d} c}{{\rm d} x} (x) \frac{{\rm d} \varphi }{{\rm d} x} (x) {\rm d}x - \frac{\lambda_-d_-}{\theta_-}  \psi (0) + \frac{\lambda_+d_+}{\theta_+}\psi (1).
$$
This equality can be rewritten as
$$
\overline{a}\left[ (u,v) \ , (\varphi ,\psi) \right] = \overline{l}\left[ (\varphi, \psi)\right] ,
$$
where $\overline{a}\left[ (u,v) \ , (\varphi ,\psi) \right]$ and $\overline{l}\left[ (\varphi, \psi)\right]$ are the left hand side and right hand side, respectively, of the previous
expression. According to this, we have the following definition of a weak solution (that we will call {\em solution} in the following).

\begin{definition}
 A (weak) solution of system (\ref{ndepefl})--(\ref{ndepsfl}) is a couple of functions $(u,v) \in H$ such that
\begin{equation} \label{edwsl}
 \overline{a}\left[ (u,v) \ , (\varphi ,\psi) \right] = \overline{l}\left[ (\varphi, \psi)\right] , \ \forall \ (\varphi , \psi)\in H.
\end{equation}
\end{definition}

It is easy to show that, in this case, $\overline{a}:H\times H \rightarrow \R$ is bilinear and continuous and $\overline{l}: H \rightarrow \R$
is linear and continuous. In order to be able to apply Lax-Milgram Theorem we need to show that $\overline{a}:H\times H \rightarrow \R$ is coercive, but this is not true (otherwise we would get
uniqueness of solution of (\ref{ndepefl})--(\ref{ndepsfl}), which, as we pointed out in Remark \ref{remseul}, is not true). Let us try to overcome this problem by changing $H$ by a suitable
closed subspace.

Let $\overline{H}=\{ \varphi\in H^1(0,L): \int_0^L\varphi(x){\rm d}x =0\} \times \left( H^1(0,L_1)\cap H^1(L_1+\delta ,1)\right)$. It is known that
$\{ \varphi\in H^1(0,L): \int_0^L\varphi(x){\rm d}x =0\}$ is a closed subspace of $H^1(0,1)$ and there exists a
constant $\alpha >0$ such that, for any function $\varphi$ in that set,
$$
\int_0^L r(x) \left( \frac{{\rm d} \varphi (x)}{{\rm d} x} (x) \right)^2  {\rm d}x \geq \alpha \int_0^L \varphi(x)^2 {\rm d}x.
$$

\begin{proposition} \label{prusl}
 There exists a unique couple of functions $(\overline{u},\overline{v}) \in \overline{H}$, satisfying
\begin{equation} \label{svpcb}
 \overline{a}\left[ (\overline{u},\overline{v}) \ , (\varphi ,\psi) \right] = \overline{l}\left[ (\varphi, \psi)\right] , \ \forall \ (\varphi , \psi)\in \overline{H}.
\end{equation}
\end{proposition}

\noindent
{\bf Proof:}
It is easy to show that $\overline{a}:\overline{H}\times \overline{H} \rightarrow \R$ is bilinear and continuous and $\overline{l}:
\overline{H} \rightarrow \R$
is linear and continuous. Hence, in order to be able to apply Lax-Milgram Theorem we need to show that $\overline{a}:\overline{H}\times \overline{H} \rightarrow \R$ is coercive. Now,

$$
\overline{a}\left[ (\varphi ,\psi) \ ,  (\varphi ,\psi) \right] =
$$
$$
\int_0^L r(x) \left( \frac{{\rm d} \varphi }{{\rm d} x} (x) \right)^2  {\rm d}x
+ \frac{\lambda_-}{\theta_-}\int_0^{L_1} \left( \frac{{\rm d} \psi}{{\rm d} x} (x) \right)^2  {\rm d}x
+ \frac{\lambda_+}{\theta_+}\int_{L_1+\delta}^L \left( \frac{{\rm d} \psi }{{\rm d} x} (x) \right)^2 {\rm d}x
$$
$$
+\lambda_- \int_0^{L_1} b_-(x) \left(  \psi (x)- \varphi (x)\right)^2  {\rm d}x
+\lambda_+ \int_{L_1+\delta}^L b_+(x) \left( \psi (x)-\varphi (x)\right)^2 {\rm d}x
$$
$$
\geq \frac{\underline{r}}{2}\int_0^L \left( \frac{{\rm d} \varphi }{{\rm d} x} (x) \right)^2  {\rm d}x
+\frac{\underline{r}\alpha }{2}  \int_0^L  \varphi(x)^2 {\rm d}x
+ \frac{\lambda_-}{\theta_-}\int_0^{L_1} \left( \frac{{\rm d} \psi}{{\rm d} x} (x) \right)^2  {\rm d}x
$$
$$
+ \frac{\lambda_+}{\theta_+}\int_{L_1+\delta}^L \left( \frac{{\rm d} \psi }{{\rm d} x} (x) \right)^2 {\rm d}x
+\lambda_-\underline{b}_- \int_0^{L_1} \left(  \psi (x)^2 + \varphi (x)^2 -2\varphi (x) \psi (x) \right)  {\rm d}x
$$
$$
+\lambda_+ \underline{b}_+\int_{L_1+\delta}^L  \left(  \psi (x)^2 + \varphi (x)^2 -2\varphi (x) \psi (x) \right)  {\rm d}x ,
$$
with ${\displaystyle \underline{r}=\inf_{x\in (0,1)}}\{r(x)\} >0$ and  ${\displaystyle \underline{b}_\pm=\inf_{x\in (0,1)}}\{b_\pm (x)\} >0$.
Using Young inequality, for any $\beta_- >0$ and $\beta_+ >0$ we get
$$
\left(  \psi (x)^2 + \varphi (x)^2 -2\varphi (x) \psi (x) \right) \geq \left(  \psi (x)^2 -\beta_\pm \psi (x)^2+ \varphi (x)^2 -\frac{1}{\beta_\pm}\varphi (x)^2 \right) .
$$
Then, if we take $\beta_-$ and $\beta_+$ defined by
$$
\frac{\underline{r}\alpha}{6} + \lambda_\pm \underline{b}_\pm = \lambda_\pm\underline{b}_\pm \frac{1}{\beta_\pm} ,
$$
we have that $\beta_- \in (0,1)$, $\beta_+ \in (0,1)$ and

$$
\overline{a}\left[ (\varphi ,\psi) \ ,  (\varphi ,\psi) \right]
\geq \frac{\underline{r}}{6}\int_0^L \left( \frac{{\rm d} \varphi }{{\rm d} x} (x) \right)^2  {\rm d}x
+ \frac{\underline{r}\alpha }{2}  \int_0^L  \varphi(x)^2 {\rm d}x
$$
$$
+ \frac{\lambda_-}{\theta_-}\int_0^{L_1} \left( \frac{{\rm d} \psi}{{\rm d} x} (x) \right)^2  {\rm d}x
+ \frac{\lambda_+}{\theta_+}\int_{L_1+\delta}^L \left( \frac{{\rm d} \psi }{{\rm d} x} (x) \right)^2 {\rm d}x
$$
$$
+\lambda_- \underline{b}_- (1-\beta_-)\int_0^{L_1} \psi (x)^2  {\rm d}x
+\lambda_+ \underline{b}_+ (1-\beta_+)\int_{L_1+\delta}^L \psi (x)^2 {\rm d}x ,
$$
$$
\geq \min \left\{ \frac{\underline{r}}{6}, \frac{\underline{r}\alpha}{2}, \frac{\lambda_-}{\theta_-}, \frac{\lambda_+}{\theta_+}, \lambda_- \underline{b}_- (1-\beta_-), \lambda_+ \underline{b}_+(1-\beta_+) \right\}
|| (\varphi , \psi)\|^2_{\overline{H}},
$$
which proves the coercivity property and finishes the proof.

\begin{theorem} \label{tmsol} If $(\overline{u} ,\overline{v})\in \overline{H}$ is the solution of (\ref{svpcb}), then $(\overline{u} + \gamma ,\overline{v} + \gamma ) \in H$ is a solution of
(\ref{ndepefl})--(\ref{ndepsfl}), for any $\gamma\in \R$.
\end{theorem}

\noindent
{\bf Proof:}
We have to prove that  $(u,v)=(\overline{u} + \gamma ,\overline{v} + \gamma )$ satisfies (\ref{edwsl}) for any $\gamma\in \R$. Now, for any $(\varphi , \psi) \in H$,
$$
(\overline{\varphi}, \overline{\psi}) =\left( \varphi - \frac{1}{L}\int_0^L \varphi (x){\rm d}x , \psi  - \frac{1}{L}\int_0^L \varphi (x){\rm d}x \right) \in \overline{H}.
$$
Then, using (\ref{svpcb}) and the fact that $\frac{\lambda_- d_-}{\theta_-}=\frac{\lambda_+ d_+}{\theta_+}=\frac{I(t)}{A}$,
\begin{eqnarray*}
& & \hspace*{-2cm} \overline{a}[(\overline{u} + \gamma ,\overline{v} + \gamma ), (\varphi , \psi)] \\
& = & \overline{a}\left[ (\overline{u} ,\overline{v} ),\left( \varphi - \frac{1}{L}\int_0^L \varphi (x){\rm d}x , \psi - \frac{1}{L}\int_0^L \varphi (x){\rm d}x \right)\right] \\
& = & \overline{l} \left( \varphi - \frac{1}{L}\int_0^L \varphi (x){\rm d}x , \psi - \frac{1}{L}\int_0^L \varphi (x){\rm d}x \right) \\
& = & \overline{l} \left( \varphi , \psi \right) + \left( \frac{\lambda_-d_-}{\theta_-}-\frac{\lambda_+d_+}{\theta_+}\right) \frac{1}{L}\int_0^L \varphi (x){\rm d}x = \overline{l} \left( \varphi , \psi \right),
\end{eqnarray*}
which concludes the proof.

\begin{corollary}
There exists a unique solution of (\ref{ndepefl})--(\ref{cfpl}).
\end{corollary}

\section{Conclusions}

A fully mathematical model for the simulation of a Lithium-ion battery is presented, based on a macro-homogeneous approach developed by Newman (see \cite{1973Ne})
It includes a system of boundary value problems for the conservation of Lithium and conservation of charge
in the solid and electrolyte phases, together with an initial value problem accounting for the conservation of energy.
These model can be very helpful for the design and optimization of new batteries and also for the real time control of its performance.

We point out that, over the last years, several authors have used similar models with an incorrect boundary condition. In this article it is showed, from a physical and a mathematical point of view,
why such a boundary condition is incorrect. To show that, we prove that the system of equations including that boundary conditions does not have any solution. Then, we deduce the correct boundary condition and we
show some results regarding the uniquesness and existence of solution of the corresponding
problem.

\section*{Acknowledgment}

This work was carried out thanks to the financial support of  the Spanish ``Ministry of Economy and Competitiveness'' under project MTM2011-22658; the research group MOMAT (Ref. 910480) supported by ``Banco Santander'' and ``Universidad Complutense de Madrid''; and the ``Junta de Andaluc\'ia'' and the European Regional Development Fund through project P12-TIC301. This publication was also supported by OCIAM (University of Oxford), where the author was collaborating as an academic visitor.

\end{document}